\documentclass[12pt]{article}
\usepackage{amsfonts,latexsym}

\title{Birational rigidity of a three-dimensional double cone
\thanks{this work was carried out with the finanical support of the Russian
Foundation for Basic Research (grant no. 96-01-00820) and a grant for
Support of Leading Scientific Schools (no. 96-15-96146)}}
\author{Mikhail Grinenko}
\date{}

\newtheorem{theorem}{Theorem}[section]
\newtheorem{proposition}[theorem]{Proposition}
\newtheorem{lemma}[theorem]{Lemma}
\newtheorem{corollary}[theorem]{Corollary}

\newtheorem{definition}[theorem]{Definition}

\makeatletter

\@addtoreset{equation}{section}
\newcommand{\l@abcd}[2]{\hbox to\textwidth{#1\dotfill #2}}

\makeatother

\newcommand*{\mybegintheorem}[1]{\begin{trivlist}\it%
      \item[\hspace{\labelsep}{\bf #1}]}
\newcommand*{\myendtheorem}{\end{trivlist}}
\newenvironment*{theorem*}{\mybegintheorem{Theorem.}}{\myendtheorem}
\newenvironment*{proposition*}{\mybegintheorem{Proposition.}}{\myendtheorem}
\newenvironment*{corollary*}{\mybegintheorem{Corollary.}}{\myendtheorem}
\newenvironment*{definition*}{\mybegintheorem{Definition.}}{\myendtheorem}

\renewcommand{\phi}{\varphi}
\renewcommand{\epsilon}{\varepsilon}

\newcommand{\lra}{\longrightarrow}

\newcommand{\PQ}{{{\mathbb P}^4}}
\newcommand{\PT}{{{\mathbb P}^3}}
\newcommand{\PTw}{{{\mathbb P}^2}}
\newcommand{\POn}{{{\mathbb P}^1}}
\newcommand{\ZA}{{\mathbb Z}}
\newcommand{\QA}{{\mathbb Q}}
\newcommand{\RA}{{\mathbb R}}

\newcommand{\freemlt}{\mathop{\rm *}}
\newcommand{\codim}{\mathop{\rm codim}\nolimits}

\newcommand{\eqdef}{\stackrel{\rm def}{=}}
\newcommand{\ord}{\mathop{\rm ord}\nolimits}
\newcommand{\supp}{\mathop{\rm supp}\nolimits}
\newcommand{\mult}{\mathop{\rm mult}\nolimits}

\newcommand{\eps}{\varepsilon}

\begin{document}

\maketitle

\paragraph{Abstract.} It is proved that a three-dimensional double cone is a
birationally rigid variety. We also compute the group of birational
automorphisms of such a variety. This work is based on the method
of "untwisting" maximal singularities of linear system.

The author would like to express his gratitude to V.A.Iskovskikh and
A.V.Pukhlikov for very useful discussions and their attention to this work.

\section{Description of a double cone}

Let $Q\subset\PQ$ be a quadratic cone, $O\in Q$ its vertex.
Suppose $R$ is a quartic in $\PQ$ such that $R_Q=R\cap Q$ is a smooth
divisor in $Q$, $O\notin R_Q$. The cone contains two pencils
${\cal P}_1$ and ${\cal P}_2$ of planes with general members $P_1$ and $P_2$.
In this situation we often write ${\cal P}_1=|P_1|$ and ${\cal P}_2=|P_2|$.

The double cover $\pi:X\stackrel{2:1}\lra Q$ branched along $R_Q$ is
called the double cone.

It is easy to see that $X$ has exactly two singularities that are ordinary
double points (both over the vertex of the cone). Further, we have two
pencils ${\cal F}_1=|F_1|$ and ${\cal F}_2=|F_2|$ of Del Pezzo surfaces
of degree 2 on X. These pencils are inverse images of ${\cal P}_1$ and
${\cal P}_2$. Indeed, the restriction of $\pi$ to $F_1\in {\cal F}_1$
gives us a double cover of $P_1\in {\cal P}_1$ branched along a smooth
quartic curve, so $F_1$ is a Del Pezzo surface with $K_{F_1}^2=2$.

There exist smooth models $\phi_U:U\to X$ and $\phi_V:V\to X$ of the
double cone, where $\phi_U$ and $\phi_V$ are birational morphisms and
isomorphisms in codimension 1 (small contractions). To construct $U$ and
$V$, we can blow up the vertex of the cone first. The exceptional divisor
is isomorphic to $\POn\times\POn$ and can be contracted along either of its
rulings onto a line. This gives us two smooth models of the cone. It remains
to take double covers of them.

The strict transform of one of the pencils, say, ${\cal F}_1$, define a
structure of Del Pezzo fibering on $V$, so we have the morphism
$\pi_V:V\to\POn$. There exist two lines $s_1$ and $s_2$ on $V$ that are
sections of $\pi_V$ and lie over the double points of $X$. The same is true
for $U$ and ${\cal F}_2$. Finally, we can make a flop in center
$s_1\cup s_2$ to get $U$ from $V$.

In the sequel, ${\cal F}_1$, $F_1$, ${\cal F}_2$ and $F_2$ will
denote objects on $X$, $V$ and $U$ simultaneously.

\paragraph{Curves on $X$, $U$ and $V$.}
It is clear that every line on $Q$ is an intersection of two elements of
the pencils ${\cal P}_1$ and ${\cal P}_2$.

\begin{definition}
A curve $C\subset X$ is called a line on $X$ if $\pi|_C$ is an isomorphism
onto a line on the cone. The lines on $V$ are the curves $s_1$, $s_2$ and the
inverse images of the lines on $X$ (Analogously, on $U$).
\end{definition}
Let $l$ be a line on $X$. Then, either $\pi(l)\subset R_Q$ or
$\pi(l)$ is bitangent to $R_Q$.

\begin{definition}
A curve $C\subset X$ is called vertical with respect to $\pi|_V$
(resp. $\pi|_U$) if $C\subset F_1$ for some $F_1\in{\cal F}_1$
(resp. $C\subset F_2$ for some $F_2\in{\cal F}_2$). $C$ is said to be
double vertical if $C\subset F_1\cap F_2$.
\end{definition}

\paragraph{Fibers of the morphisms $\pi|_V$ and $\pi_U$.}
Everywhere below we assume that the following condition holds:
\begin{itemize}
\item[{\bf(A)}]  elements of ${\cal F}_1$ or ${\cal F}_2$ are smooth
    along any double vertical line
\end{itemize}
This condition always holds for a general double cone because of the finitely
many nonsmooth elements of the pencils ${\cal F}_1$ and ${\cal F}_2$.

\paragraph{The divisor class group of $X$.}
Suppose $S$ is a small resolution of the singularity of the cone such that
$\pi:V\to S$ is a double cover. Let $P_1$ and $P_2$ be the classes in $Cl(S)$
of the strict transforms of the respective planes in $Q$. By $s$ denote
the exceptional line on $S$. We obviously have $P_1\circ s=1$ and
$P_2\circ s=-1$. Then,
$$
       Pic(S)=\ZA[P_1]\oplus\ZA[P_2]
$$
and $K_S\sim -3H$, where $H=P_1+P_2$ is a lift of a hyperplane section of
the cone.
\begin{proposition}
  $Pic(V)\simeq \ZA\oplus\ZA$.
\end{proposition}
{\bf Proof.} It is obvious that $(-K_S)$ and $(-K_V)\sim\pi^*H$ are
{\it big-} and {\it nef-} divisors. Using the Serre duality and the
well-known vanishing theorem for {\it big-} and {\it nef-} divisors, we have
$$
   h^{1,0}(S)=h^{2,0}(S)=h^{1,0}(V)=h^{2,0}(V)=0.
$$
Combining that with cohomologies of the exponential sequence, we get
$$
     Pic(S)\simeq H^2(S,\ZA)
$$
and
$$
  Pic(V)\simeq H^2(V,\ZA).
$$
To prove the proposition, we only need to show that $h^{1,1}(V)=h^{1,1}(S)$.

Suppose ${\cal L}={\cal O}_S(2H)$. We may assume that the morphism
$\pi:V\to S$ is branched along the divisor $R\in\Gamma(S,{\cal L}^2)$.
Then, let
$$
    L={\bf Spec}\,(Sym\, {\cal L}^{-1})
$$
be the total space of the line bundle ${\cal L}$ with the projection
$p:L\to S$. If $f=0$ is a local equation of $R$, then
${t^2=f}$ define $V\subset L$ and $\pi=p|_V$.

We have the following exact sequence on V:
$$
   0\lra\pi^*{\cal O}_S(-R)\lra\Omega_L|_V\lra\Omega_V\lra 0.
$$
Note that $\pi^*{\cal O}(R)$ is {\it big} and {\it nef} on $V$. Therefore,
we have isomorphisms
$$
   H^1(V,\Omega_V)\simeq H^1(V,\Omega_L|_V)\simeq H^1(S,\pi_*\Omega_L|_V)
$$
The last follows from vanishing $R^i\pi_*=0$ for any $i>0$.

Further, there is an exact triple of sheaves on $L$:
$$
   0\lra p^*\Omega_S\lra\Omega_L\lra p^*{\cal L}^{-1}\lra 0
$$
We can restrict it to $V$, and then get the direct image on $S$:
$$
   0\lra\Omega_S\oplus\Omega_S(-2H)\lra\pi_*\Omega_L|_V\lra%
      {\cal O}_S(-2H)\oplus{\cal O}_S(-4H)\lra 0
$$
Using the vanishing theorem again, we obtain
$$
  H^1(S,\pi_*\Omega_L|_V)\simeq H^1(S,\Omega_S)\oplus H^1(S,\Omega_S(-2H))
$$
It remains to show that $H^1(S,\Omega_S(-2H))=0$. There is an exact sequence
on $F\in{\cal P}_2$:
$$
  0\lra\Omega_S(-2H)\lra\Omega_S(-2H+F)\lra%
         \Omega_S(-2H+F)\otimes{\cal O}_F\lra 0
$$
Notice that $H^1(S,\Omega_S(-2H+F))=0$, since the divisor
$2H-F\sim P_2+2P_1$ is ample on $S$.
But $H^0(F,\Omega_S(-2H+F)\otimes{\cal O}_F)$ vanishes, too. Indeed,
that follows from the sequence
$$
  0\lra{\cal N}^{-1}_{F|S}\otimes{\cal O}_S(-2H+F)\lra\Omega_S(-2H+F)|_F\lra%
     \Omega_F(-2H+F)\lra 0
$$
Thus, $H^1(S,\Omega_S(-2H))=0$. The proposition is proved.

\paragraph{} We can see now that
$$
   Cl(X)=\ZA[F_1]\oplus\ZA[F_2],
$$
and the same is true for $U$.

We recall that there exist the lines $s_1$ and $s_2$ on $V$, which
are sections of the fibering $\pi_V:V\to \POn$. Denote a vertical line
class by $f$. Then the 1-cycle group
$$
      A^2_\QA\eqdef A^2\otimes\QA
$$
is generated by $s_1$ and $f$ (or $s_2$ and $f$, this is the same), i.e.,
$$
      A^2_\QA=\QA s_1\oplus\QA f
$$

\section{Description of birational automorphisms}
\paragraph{Birational automorphisms associated with curves on $X$.}
Consider the smooth model $\phi_V:V\to X$ of the double cone.
Let $l\subset V$ be a section of the fibering $\pi_V:V\to\POn$ that is not
$s_1$ or $s_2$ and is not contained in the ramification divisor. Then there
is a curve $l^*\subset V$ conjugate to $l$ with respect to the involution
transposing the sheets of the cover $\pi:V\to Q$. Note that $l^*$ is a
section too. Let $\psi:\tilde V\to V$ be the composition of a blow-up of
the curve $l$ and then of the proper inverse image of $l^*$. We denote by
$E$ and $E^*$ the corresponding exceptional divisors on $\tilde V$.

We now describe the birational automorphisms $\tau_l$ associated with the
section $l$. By definition, the restriction of $\tau_l$ to a general
fiber $T\in |F_1|$ is Bertini's involution. In detail, we can blow up
the points $A=l\cap T$ and $A^*=l^*\cap T$ to see $T$ as a fibering on
elliptic curves $\tilde T$. The corresponding exceptional curves $e$ and
$e^*$ are the sections of this fibering. The fiber reflection with respect
to $e^*$ gives us an automorphism $\tilde\tau_T$ that lowers to $T$ as
Bertini's involution $\tau_T$ (see \cite{Isk4}). Thus,
$$
     \tau_l|_T\eqdef\tau_T.
$$

Consider $\tilde T$ as a surface in $\tilde V$, so that
$e=E|_{\tilde T}$ and $e^*=E^*|_{\tilde T}$. Suppose
$h=-K_{\tilde V}|_{\tilde T}$. It can easily be checked that
$$
\left\{
\begin{array}{l}
\tilde\tau^*_T(e)=-e+2e^*+2h\\
\tilde\tau^*_T(e^*)=e^*\\
\tilde\tau^*_T(h)=h
\end{array}
\right.
$$
Let we have a linear system $|D|\subset |-nK_V+mF_1|$ on $V$ such that
$n>0$ and $m\ge 0$. Suppose $\nu=\mult_l|D|$ and $\nu^*=\mult_{l^*}|D|$;
then
$$
\tilde D\sim\psi^*(-nK_V+mF_1)-\nu E-\nu^* E^*,
$$
where $\tilde D$ is the strict transform of $D$ on $\tilde V$. Further,
we can write
$$
  \tilde\tau^*_l(\tilde D)\sim\psi^*(-aK_V+bF_1)-cE-dE^*.
$$
To compute unknown coefficients, it is useful to restrict all to $\tilde T$.
Suppose $h'=h+e+e^*=\psi^*(-K_V)|_{\tilde T}$; we have
$$
\begin{array}{l}
\tilde D|_{\tilde T}=nh'-\nu e-\nu^*e^*\\
\tilde\tau^*_l(\tilde D)|_{\tilde T}=ah'-ce-de^*
\end{array}
$$
Note that
$$
  \tilde\tau^*_l(\tilde D)|_{\tilde T}=%
     \tilde\tau^*_{\tilde T}(\tilde D|_{\tilde T}).
$$
Using the action of $\tilde\tau^*_T$ on generators of $Pic(\tilde T)$,
we get
$$
\tilde\tau^*_{\tilde T}(\tilde D|_{\tilde T})=%
(3n-2\nu)h'-(4n-3\nu)e^*-\nu^*e^*
$$
So, $a=3n-2\nu$, $c=4n-3\nu$ и $d=\nu^*$. We will not use the
coefficient $b$.

Thus,
\begin{equation}
\label{tau_l}
  \tau_l^*(D)\in |(3n-2\nu)(-K_V)+\ldots F_1-(4n-3\nu)l-\nu^*l^*|
\end{equation}

In the same way, we define birational automorphisms associated with
sections of the fibering $\pi_U:U\to\POn$.

\paragraph{Birational automorphisms associated with double points on $X$.}
There are another two birational automorphisms. Let $E_1$ and $E_2$ be
exceptional divisors after blowing up both singular points
$\phi:\tilde V\to X$. The composition of $\phi$, $\pi$, and a projection
$Q\to\POn\times\POn$ from the vertex of the cone realize $\tilde V$ as
a fibering on elliptic curves $\tilde V\to\POn\times\POn$ such that
$E_1$ and $E_2$ are its sections. Suppose $\tilde\tau_1$ and $\tilde\tau_2$
are fiber reflections with respect to the divisors $E_1$ and $E_2$.
The birational morphisms are
$$
\begin{array}{c}
\tau_1=\phi\circ\tilde\tau_1\circ\phi^{-1}\\
\tau_2=\phi\circ\tilde\tau_2\circ\phi^{-1}
\end{array}
$$
It is easy to see that
$$
   Pic(\tilde V)=\ZA[\tilde F_1]\oplus\ZA[\tilde F_2]%
       \oplus\ZA[E_1]\oplus\ZA[E_2],
$$
where $\tilde F_1$ and $\tilde F_2$ are the strict transforms of
$F_1$ and $F_2$. Note that elements of the pencils $|\tilde F_1|$ and
$|\tilde F_2|$ are invariant with respect to $\tilde\tau_1$ and
$\tilde\tau_2$.

Let $\tilde S\in |\tilde F_1|$ be a general element. By $g$ denote a fiber
of a natural projection  $\tilde S\to\POn$. Note that $e_1=E_1\cap\tilde S$
and $e_2=E_2\cap\tilde S$ are sections of this projection. Obviously,
$e_1^2=e_2^2=-1$, $g^2=0$, $e_1\circ e_2=0$, and $e_1\circ g=e_2\circ g=1$.
Suppose $\tau_{1_{\tilde S}}\eqdef\tilde\tau_1|_{\tilde S}$.
It is easy to check that
$$
\left\{
\begin{array}{l}
\tau_{1_{\tilde S}}^*(e_1)=-e_1+2e_2+2g\\
\tau_{1_{\tilde S}}^*(e_2)=e_2\\
\tau_{1_{\tilde S}}^*(g)=g
\end{array}
\right.
$$
These formulae are also true if $\tilde S\in |\tilde F_2|$.

Now, consider the variety $V$. It is clear that $\tilde V$ is a blow-up
$\psi:\tilde V\to V$ of the lines $s_1$ and $s_2$ on $V$. Let we have a
linear system
$$
    |D|\subset |-nK_V+mF_1-\nu_1 s_1-\nu_2 s_2|
$$
on $V$, where $n>0$, $m\ge 0$, $\nu_1=\mult_{s_1}|D|$, and
$\nu_2=\mult_{s_2}|D|$. Denote a lift of $\tau_1$ to $V$ by the same letter.
The strict transform of the linear system $|D|$ on $\tilde V$ is
$$
    |\tilde D|\subset |\psi^*(-nK_V+mF_1)-\nu_1 E_1-\nu_2 E_2|.
$$
Suppose
$$
    \tau_1^*(\tilde D) = |\psi^*(-aK_V+bF_1)-c E_1-d E_2|.
$$
Restricting $\tilde D$ to $\tilde S\in |\tilde F_1|$ and then to
$\tilde S\in |\tilde F_2|$ give us $a=3n-2\nu_1$, $b=m$, $c=4n-3\nu_1$,
and $d=\nu_2$. After lowering to $V$ we get
\begin{equation}
\label{tau_1}
  \tau_1^*(D)\in |(3n-2\nu)(-K_V)+mF_1-(4n-3\nu)s_1-\nu_2s_2|
\end{equation}
The automorphism $\tau_2$ is analoguosly described. In fact, $\tau_1$
and $\tau_2$ are conjugate in the group of birational automorphisms by
the regular involution $\delta: X\to X$ transposing the sheets of the cover
$\pi: X\to Q$.

\section{Formulation of the main result}
\begin{definition}
A pair $(W,{\cal Y})$, where $W$ is nonsingular in codimension 1 progective
variety and $\cal Y$ is a linear system on $W$ is called test if the
following conditions hold:
\begin{itemize}
\item[(i)] ${\cal Y}=|Y|$ is free from fixed components;
\item[(ii)] there is a number $\alpha=\alpha(W,{\cal Y})\in\RA_+$ such that
for every rational number $\beta > \alpha$ the linear system
$$
   |m(Y+\beta K_W)|
$$
is empty for any $m\in\ZA_+$, $m\beta\in\ZA$.
\end{itemize}
\end{definition}
The indicated number $\alpha(W,{\cal Y})$ is called {\it the adjunction
treshold} of the linear system $\cal Y$.
\begin{definition}
A variety $S$ is said to be birationally rigid if for any test pair
$(W,{\cal Y})$ and for any birational map $\chi: S-\!\! -\!\!\to W$
there is a birational automorphism of $S$ such that
$$
         \alpha(S,{\cal T})\le\alpha(W,{\cal Y}),
$$
where ${\cal T}=(\chi\circ\chi')^{-1}{\cal Y}$.
\end{definition}

\begin{theorem}
\label{main_theorem}
Under the conditions of section 1, the double cone is birationally rigid.
\end{theorem}

\paragraph{} By $B(X)$ denote a group free generated by all
birational automorphisms of kind $\tau_l$, $\tau_1$ or $\tau_2$.
Further, let $I_X\subset Bir(X)$ be a group generated by the regular
involution $\delta$ and the birational automorphism $\tau_1$.
By $V_{\eta}$ denote the fiber over the common point of the fibering
$\pi_V:V\to\POn$, and by $U_{\eta}$ the same for $\pi_U:U\to\POn$.
Obviously, there is the natural embedding of $I_X$ in the groups of
birational automorphisms $Bir(V_{\eta})$ and $Bir(U_{\eta})$.

\begin{corollary}
\label{main_corollary}
(i) The double cone is not rational and not isomorphic to a conic budle.

\noindent (ii) The group of birational automorphisms $Bir(X)$ is the
semi-direct product of $B(X)$ and the group of regular automorphisms
$Aut(X)$ of the double cone, i.e., the following exact triple holds:
\begin{equation}
\label{sl1}
      0\lra B(X)\lra Bir(X)\lra Aut(X)\lra 0
\end{equation}
Moreover, in general case, $Bir(X)$ is the free product of $Bir(V_{\eta})$
and $Bir(U_{\eta})$ with the common subgroup $I_X$:
\begin{equation}
\label{sl2}
       Bir(X)\simeq Bir(V_{\eta})\freemlt_{I_X} Bir(U_{\eta})
\end{equation}
\end{corollary}

\section{ Maximal singularities}
\label{MAXSing}
Let we have a test pair $(W,{\cal Y})$ and a birational map
$$
         \chi: X-\!\! -\!\!\to W.
$$
The idea of proving the theorem is the following. Assume
$\alpha(X,|D|)>\alpha(W,{\cal Y})$, where $|D|=\chi^{-1}{\cal Y}$.
Choose the smooth model of $X$ (either $V$ or $U$) such that the strict
transform of the linear system $|D|$ is
$$
   |D|\subset |-\alpha(X,|D|)K + \{\mbox{ fibers }\}|.
$$
For example, let it be $V$. Then
\begin{equation}
\label{Dform}
  |D|\subset |-nK_V+mF_1|,
\end{equation}
where $n=\alpha(X,|D|)>0$ and $m\ge 0$.

Then using $\tau_l$, $\tau_1$ или $\tau_2$, we show that the linear system
$|D|$ can be reduced to to form
$$
  |D'|\subset |-n'K_V+m'F_1|,
$$
where $n'<n$, so $\alpha(X,|D'|)<\alpha(X,|D|)$. Then we lower to $X$ and
choose the smooth model again, but with respect to $|D'|$, and so further.
After some steps we will obtain the statement of the theorem.

The first fact we need is an existence of maximal singulatities of linear
system under conditions as above (see \cite{Pukh2},\cite{Pukh3}).
Suppose $\chi_1$ is a lift of $\chi$ to $V$. Let $|D|=\chi_1^{-1}{\cal Y}$
have the form (\ref{Dform}),
$n=\alpha(V,|D|)=\alpha(X,|D|)>\alpha(W,{\cal Y})\ge 0$.
Then there exist a smooth quasi-projective variety $\tilde V$ and
a birational map $\phi:\tilde V\to V$ such that
\begin{equation}
\label{empty}
 \emptyset=|\phi^{-1}(D)+nK_{\tilde V}|=%
    |m\phi^*(F_1)+\sum_{T\in{\cal L}}(n\delta(T)-\nu_T(|D|)T)|
\end{equation}
where $\cal L$ is a set of Weil divisors on $V$ that are exceptional for
$\phi$,
$$
 K_W=\phi^*(K_V)+\sum_{T\in{\cal L}}\delta(T)T
$$
($\delta(T)$ is so-called discrepancy), and $\nu_T(|D|)=\ord_T(\phi^*(|D|))$.

Suppose $e_T=\nu_T(|D|)-\delta(T)$. We say that the linear system $|D|$
has {\it a maximal singularity} with a center $B_T=\phi(T)$, if $e_T>0$.
We also say that the valuation $\nu_T$ realize this maximal singularity.

It follows from (\ref{empty}) that the set of maximal singularities
$$
   {\cal M}=\{T\in{\cal L}: e_T>0\}
$$
is not empty.

The most difficult part of this paper is to prove the following fact:
\begin{proposition}
\label{difficult_prop}
Centers of maximal singularities cannot be points only.
\end{proposition}

\paragraph{Proof.} The proof is by {\it reductio ad absurdum.}

We begin with some definitions.
\begin{definition}
Suppose $W$ is a smooth projective variety, $C$ is a curve, and $G$ is a
divisor on $W$ such that $\supp C\not\subset\supp G$. A degree of $C$ with
respect to $G$ is $\deg_GC=G\circ C$.
\end{definition}

\begin{lemma}
Suppose a point $x\in C\cap G$; then $\mult_xG\cdot\mult_xC\le\deg_GC$.
\end{lemma}
The proof is trivial.

We recall that $V$ has the fiber structure $\pi_V:V\to\POn$.
\begin{definition}
A curve $C$ is called horizontal (with respect to $\pi_V$)
if any component of $C$ covers the base (i.e., $\POn$).
\end{definition}

Now, let $D_1$ and $D_2$ be general elements of $|D|$. Suppose
$$
     D_1\circ D_2=\alpha_1s_1+\alpha_2s_2+Z^h+\sum_{t\in\POn}Z_t^v,
$$
where $Z^h$ is a horizontal 1-cycle, $s_1, s_2\not\subset\supp Z^h$,
and $Z_t^v$ is a vertical 1-cycle in the fiber over a point $t\in\POn$.
Note that all cycles in this decomposition are effective, since $|D|$
has no fixed components.

Choose an element $G\in |F_2|$ that does not contain components of $Z^h$
or $Z_t^v$. Then
$$
  D_1\circ D_2\circ G=2n^2+4mn\ge -\alpha_1-\alpha_2+%
                    \sum_{t\in\POn}deg_GZ_t^v.
$$
By $S_t$ denote a fiber over a point $t\in\POn$. From the condition
(\ref{empty}) it follows that
$$
    m < \sum_{t\in\POn} \max_{\{T:\phi(T)\in S_t\}}\frac{e_T}{\nu_T(S_t)}
$$
(see \cite{Pukh4}, proposition 2.1).
So we have
$$
   \sum_{t\in\POn}deg_GZ_t^v < 2n^2+\alpha_1+\alpha_2+%
        4n\sum_{t\in\POn} \max_{\{T:\phi(T)\in S_t\}}\frac{e_T}{\nu_T(S_t)}
$$
Let $p$ be the number of fibers $\pi_V$ with maximal singularities of $|D|$.
Then there exists a maximal singularity $T$, $\phi(T)=t$, such that the
following inequality holds:
$$
   \nu_T(S_t)\deg_GZ_t^v < \frac{2n^2+\alpha_1+\alpha_2}{p}\nu_T(S_t)+4ne_T
$$
This is the supermaximal singularity in notation of paper \cite{Pukh4}.

Suppose $B_0=\phi(T)$, $S=S_t$, $e=e_T$, $\nu=\nu_T$, $\delta=\delta_T$,
and take a resolution of the maximal singularuty in the usual way
(\cite{Pukh3}). The chain
$$
 V_N\stackrel{\phi_{N,N-1}}{\longrightarrow}V_{N-1}%
 \stackrel{\phi_{N-1,N-2}}{\longrightarrow}\ldots%
 \stackrel{\phi_{i+1,i}}{\longrightarrow}V_i%
 \stackrel{\phi_{i,i-1}}{\longrightarrow}\ldots%
 \stackrel{\phi_{2,1}}{\longrightarrow}V_1%
 \stackrel{\phi_{1,0}}{\longrightarrow}V
$$
is the composition of blow-ups with centers $B_i$ and with exceptional
divisors $E_i\subset V_i$ such that for every $i$ $\phi_{i,i-1}(B_i)=B_{i-1}$
and a triple $(V_N,E_N,\phi_{N,0})$ realize the valuation $\nu$, where
$\phi_{N,0}=\phi_{N,N-1}\circ\ldots\circ\phi_{1,0}$. We may assume that
$B_0,\ldots B_{L-1}$ are points, and $B_L,\ldots B_{N-1}$ are curves,
$1\le L<N$.

Let $\Gamma$ be the directed graph of the singularity: vertices $i$
and $j$ are joined by an arrow (from $i$ to $j$) if $i>j$ and
$B_{i-1}\subset E_j^{i-1}$ (upper indexes indicate strict transforms on the
corresponding floors of the chain of blow-ups). We denote by $r_i$
the number of different paths from the vertex $N$ to a vertex $i$.
It is not very difficult to check this relation:
\begin{equation}
\label{r_i}
      r_i = \sum_{j\to i} r_j
\end{equation}
Suppose
$$
\begin{array}{l}
       \nu_i=\mult_{B_{i-1}}|D|^{i-1}, \\
       \delta_i=\codim B_{i-1}-1,
\end{array}
$$
and $\Sigma_0=\sum_{i=1}^Lr_i$, $\Sigma_1=\sum_{i=L+1}^Nr_i$.
So,
$$
\begin{array}{l}
   \nu=\sum_{i=1}^Nr_i\nu_i, \\
   \delta=2\Sigma_0+\Sigma_1.
\end{array}
$$
We obtain so-called Noeter-Fano-Iskovskikh's inequality:
$$
    e=\nu-n\delta=\sum_{i=1}^Nr_i\nu_i-2n\Sigma_0-n\Sigma_1 > 0.
$$
Since $\nu_1\ge\nu_2\ge\ldots\ge\nu_N$, then $\nu_1 > n$. From here it
follows that $S$ is smooth at the point $B_0$. Indeed, otherwise
$\mult_{B_0}S=2$, and then $S$ is a double cover of a plane branched
along a quartic with a singular point at $B_0$. Take a general line on the
plane through this point, and let $C\subset S$ be the inverse image of $l$.
We obtain a contradiction:
$$
     2n = D\circ C\ge \mult_{B_0}C\cdot\mult_{B_0}|D| = 2\nu_1 > 2n.
$$

Further, let us fix the number $L'=\max\{i\le L: B_{i-1}\in S^{i-1}\}$,
then
$$
     \Sigma_0'\eqdef r_1+\ldots+r_{L'}\le \nu(S)
$$
Suppose $m_{0,i}^v=\mult_{B_{i-1}}(Z^v)^i$ and
$m_{0,i}^h=\mult_{B_{i-1}}(Z^h)^i$. It is obvious that
$m_{0,i+1}^h\le m_{0,i}^h$ and $m_{0,i+1}^v\le m_{0,i}^v$ for any $i$.
Note also that
$$
    m_{0,1}^h\le 2n^2-\alpha_1-\alpha_2.
$$
Now the supermaximal singularity condition is
$$
    \Sigma_0'm_{0,1}^v <
    \Sigma_0'\left(\frac{2n^2+\alpha_1+\alpha_2}p\right)+4ne
$$
Assume that the point $B_0$ is not contained in $s_1$ or $s_2$. Then we
can deduce from the Noether-Fano-Iskovskikh inequality in the usual way that
\begin{equation}
\label{nermult}
     \Sigma_0m_{0,1}^h+\Sigma_0'm_{0,1}^v \ge \sum_{i=1}^Nr_i\nu_i^2 \ge
                   \frac{(2n\Sigma_0+\Sigma_1+e)^2}{\Sigma_0+\Sigma_1}
\end{equation}
Since $p\ge 1$ and $\Sigma_0\ge\Sigma_0'$, we get
$$
  (\Sigma_0+\Sigma_1)(4n^2\Sigma_0+4ne) >
  4n^2(\Sigma_0+\Sigma_1)\Sigma_0+n^2\Sigma_1^2+2ne(2\Sigma_0+\Sigma_1)+e^2
$$
So we have a contradiction:
$$
       (n\Sigma_1-e)^2 < 0
$$
\noindent {\bf Remark.} To this moment, we essentially followed
to arguments of paper \cite{Pukh4}.

Now suppose $B_0\in s_1$ (the case $B_0\in s_2$ is analogous). Then we have
to take into account the multiplicity along the line $s_1$ (I mean to add
$\alpha_1\Sigma_0$ to the left part of the inequality (\ref{nermult})).
After some elementary transformations we obtain
$$
  (\Sigma_1n-e)^2+\Sigma_0(\Sigma_0+\Sigma_1)%
         \left[2n^2-\left(\frac{2n^2+\alpha_1}{p}\right)%
              \cdot\frac{\Sigma_0'}{\Sigma_0}\right] < 0,
$$
hence
$$
    \frac{2n^2+\alpha_1}p > \frac{2n^2\Sigma_0}{\Sigma_0'}.
$$
Since $\Sigma_0'\le\Sigma_0$ and $\alpha_1\le 2n^2$, then $p=1$ and
\begin{equation}
\label{n1}
       \Sigma_0' > \frac12\Sigma_0.
\end{equation}
So, only one fiber of $\pi_V$ can contain maximal singularities of the
linear system $|D|$.

\section{ Infinitely near singularities}
\label{BB1}
It remains to consider infinitely near cases, i.e., $B_0\in s_1$ or
$B_0\in s_2$.

Since only one fiber of $\pi_V$ contains maximal singularities of $|D|$,
from (\ref{empty}) it follows that the {\it stronger} Noether-Fano-Iskovskikh
inequality for the maximal singularity over the point $B_0\in s_1$ holds,
i.e.,
$$
    \nu(|D|) > n\delta + m\nu(S)
$$
So,
\begin{equation}
\label{NFI}
   \sum_{i=1}^Nr_i\nu_i > 2n\Sigma_0 + n\Sigma_1 + m\Sigma_0'
\end{equation}
Take the number
$$
   q = \max \{i=1\ldots L: B_{i-1}\in s_1^{i-1}\}
$$
Suppose $q >1 $; so  $L > 1$. Then we get $\nu(S)=\Sigma_0'=r_1$, since
$s_1^1\cap S^1 = \empty$.

Assume first that $q\ge 3$; we have $B_1\in s_1^1$ and $B_2\in s_1^2$.
This means that $B_2\notin E_1^2$, hence $r_2=r_1$ from (\ref{r_i}).
Then
$$
    \frac{\Sigma_0'}{\Sigma_0} = \frac{r_1}{r_1+r_2+\ldots} \le \frac12,
$$
but this is impossible because of (\ref{n1}).

We recall the notation of the maximal singulariries method (see \cite{Pukh3}).
Let $D_1$ and $D_2$ be general elements of the linear system $|D|$, and
$$
    D_1\circ D_2 = \alpha_1s_1+\alpha_2s_2+Z_0^h+Z_0^v,
$$
where $Z_0^h$ and $Z_0^v$  are horizontal and vertical,
$s_1,s_2\not\subset\supp Z_0^h$.

Suppose
$$
    D_1^i\circ D_2^i = \alpha_1s_1^i+\alpha_2s_2^i+(Z_0^h)^i+(Z_0^v)^i+%
             Z_1^i+\ldots+Z_{i-1}^i+Z_i,
$$
where $Z_k\subset E_k$ is a curve; by $d_k$ denote either the degree of this
curve if $E_k\simeq\PTw$, or the intersection index of $Z_k$ with a fiber
of $E_k$ if $E_k$ is a ruled surface. Suppose
$m_{i,j}\eqdef\mult_{B_{j-1}}Z_i^{j-1}$.

We can decrease the coefficients $r_1,\ldots r_L$ a little, that will be
very useful for the sequel.Further we will suppose
\begin{equation}
\label{truk}
          r_i = \sum_{{j\to i}\atop{j\le L}}r_j
\end{equation}

Assume now that $q=2$. Following to paper \cite{Pukh3},\S 7, we can
write that
$$
\left\{
\begin{array}{rcl}
   \alpha_1+m_{0,1}^h+m_{0,1}^v & = & \nu_1^2+d_1\\
   \alpha_1+m_{0,2}^h+m_{1,2} & = & \nu_2^2+d_2\\
                 \cdots &&\\
   m_{0,L}^h+m_{1,L}+\ldots m_{L-1,L} & = & \nu_L^2+d_L
\end{array}
\right.
$$
The standard trick with annihiliation of the quantities $d_i$ and $m_{i,j}$
still works, despite the decrease of the coefficients $r_k$, and, using
the quadratic inequality, we get that
$$
 (r_1+r_2)\alpha_1+\sum_{i=1}^Lm_{0,i}r_i+r_1m_{0,1}^v >
           \frac{(2n\Sigma_0+n\Sigma_1+mr_1)^2}{\Sigma_0+\Sigma_1}
$$
Here the right-hand side is greater than $4n^2\Sigma_0+4mnr_1$. Taking into
account that
\begin{equation}
\label{n2}
\begin{array}{c}
       m_{0,1}^h+m_{0,1}^v\le \deg_{(-K_V)}(D_1\circ D_2) = 4n^2+4mn,\\
       m_{0,1}^h \le 2n^2-\alpha_1,
       \end{array}
\end{equation}
we get a contradiction:
$$
   (r_1+r_2)\alpha_1 > 2n^2r_1+\alpha_1r_1,
$$
i.e., $\alpha_1 > 2n^2$.

Suppose now $q=1$ and $B_1$ (i.e., $L>1$) is a point that is not
contained in $s_1$. Arguing as before, we obtain
$$
 r_1\alpha_1+\sum_{i=1}^Lm_{0,i}^hr_i+r_1m_{0,1}^v > 4n^2\Sigma_0+4mnr_1,
$$
whence, using (\ref{n2}), we have a contadiction again:
$$
  r_1\alpha_1+2n^2\sum_{i=2}^Lr_i\ge r_1\alpha_1+\sum_{i=2}^Lr_im_{0,i}^h >
      4n^2\sum_{i=2}^Lr_i,
$$
i.e., $r_1\alpha_1 > 2n^2\sum_{i=2}^Lr_i\ge 2n^2r_1$.

Two cases remain: either $B_1\in S^1$, or $B_1\subset E_1$ is a line
that is not contained in $S^1$. Indeed, the inequality (\ref{NFI})
give us $\nu_1+\nu_2 > 2n$, and, if the line $B_1\subset S^1$, we obtain
a contradiction for the strict transform of a curve $C\in |-K_S-B_0|$:
$$
   C^1\circ D^1=2n-\nu_1\ge\nu_2
$$

So, at the point $B_0$ we have an inequality for multiplicities of cycles
in $D_1\circ D_2$:
$$
  \alpha_1+m_{0,1}^h+m_{0,1}^v >
\frac{(2nr_1+n\Sigma_1+mr_1)^2}{(r_1+\Sigma_1)r_1} = 4n^2+4mn+\phi(\Sigma_1),
$$
where
$$
   \phi(\Sigma_1) = \frac{(n\Sigma_1-mr_1)^2}{(r_1+\Sigma_1)r_1}
$$
Suppose $2\theta=\nu_1+\nu_2$. From (\ref{NFI}) it follows that
$$
     \Sigma_1 > \frac{2n-\theta}{\theta-n}r_1+\frac{m}{\theta-n}r_1
$$
Hence, since $\phi'(\Sigma_1) >0$, we get
$$
  \phi(\Sigma_1) >
  \phi\left(\frac{2n-\theta}{\theta-n}r_1+\frac{m}{\theta-n}r_1\right)=
        \frac{(2n-\theta)^2}{\theta-n}(n+m)
$$
\begin{lemma}
$$
          \nu_1+\nu_2 \le \frac52n
$$
\end{lemma}
{\bf Proof.} Suppose either $\tilde B_1=B_1$ if $B_1$ is a point, or
$\tilde B_1=B_1\cap S^1$ if $B_1$ is a curve.

Let $L\in |-K_S|$ pass  through the points $B_0$ and $\tilde B_1$
(the second is infinitely near the first). Suppose $L$ is reducible,
$$
         L=L_1+L_2,
$$
where $L_1$ and $L_2$ are conjugate to each other (-1)-lines. We may assume
that only $L_1$ pass through $B_0$ and $\tilde B_1$. Recall that $S^1$ is
smooth along these lines from  condition {\bf A}.

We have (\cite{Pukh3}, lemma 6.1)
\begin{equation}
\label{locform}
    D^1\circ S^1 = (D\circ S)^1+mE_{S^1},
\end{equation}
where $E_{S^1}=E_1\cap S^1$, $m\ge 0$. Suppose $\nu=mult_{B_0}(D|_S)$,
$\tilde\nu=mult_{\tilde B_1}(D^1|_{S^1})$, and
$$
     D|_{S} =  C+k_1L_1+k_2L_2,
$$
where $C$ is an effective divisor on $S$, $L_1,L_2\not\subset\supp C$.
Then
$$
\begin{array}{rcl}
\nu & = & mult_{B_0}C+k_1, \\
\tilde\nu & = & mult_{\tilde B_1}C^1+k_1+m,
\end{array}
$$
hence
$$
   \nu+\tilde\nu = mult_{B_0}C+mult_{\tilde B_1}C^1+2k_1+m \le %
          C\circ L_1+2k_1+m
$$
Since $C\sim n(-K_S)-k_1L_1-k_2L_2$, we have $C\circ L_1=n+k_1-2k_2$. So
$$
   \nu+\tilde\nu \le n+3k_1-2k_2+m
$$
Using $C\circ L_2=n-2k_1+k_2\ge 0$, we get
$$
   \nu+\tilde\nu \le \frac52n+m
$$
It remains to take into account that $\nu=\nu_1+m$ from (\ref{locform}),
and $\tilde\nu\ge\nu_2$.

Now let $L$ be irreducible. Suppose $D|_{S} =  C+kL$.
Arguing as above, we obtain
$$
   \nu+\tilde\nu \le C\circ L+2k+m,
$$
but $C\circ L=2n-2k$. This proves the lemma.

From the lemma, $\phi(\Sigma_1) > 2n^2$, whence
$$
  \alpha_1+m_{0,1}^h+m_{0,1}^v > 6n^2+4mn
$$
On the other hand,  $m_{0,1}^h+m_{0,1}^v\le 4n^2+4mn$, and we obtain a
contradiction: $\alpha_1>2n^2$.

Proposition \ref{difficult_prop} is completely proved.

\section{Maximal curves}
\label{MC}
So, the linear system $|D|\subset |-nK_V+mF_1|$ cannot have maximal
singularities over points.
\begin{proposition}
(i) Only section of $\pi_V:V\to\POn$ that are not contained in the
ramification divisor can be centers of maximal singularities.

\noindent (ii) If $|D|\subset |-nK_V|$, then the maximal curve
is either the line $s_1$ or $s_2$.

\noindent (iii) $D$ cannot have maximal singularities in two curves
simultaneously (including an infinitely near case).
\end{proposition}
{\bf Proof.} {\it (i)} This case is given in \cite{Pukh4}, \S 4.

\noindent {\it (ii)} Let $B$ be a maximal curve of the linear system
$|D|\subset |-nK_V|$. From {\it (i)} it follows that the class of $B$
in $A^2_{\QA}$ is $s+\alpha f$, where $\alpha\ge 0$. Suppose
$\nu=mult_B|D|$ and
$$
    D_1\circ D_2 = \eps B+C,
$$
where $C$ is the effective 1-cycle, $B\not\subset\supp C$ and
$\eps\ge\nu^2$. For a general element $H\in |-K_V|$ we have
$$
   4n^2=D_1\circ D_2\circ H \ge B\circ H = \alpha\eps,
$$
whence $\alpha < 4$. It need to prove that $\alpha$ is equal to 0 or 1.

assume the converse. Since $B$ is a section that is not contained in the
ramification divisor, there exists the curve $B^*$ conjugate to $B$
with respect to the involution $\delta$ transposing the sheets of the
cover $\pi:V\to Q$. It is known also that the image of $B$ on $Q$ is
either a conic (for $\alpha=2$), or a space cubic curve (for $\alpha=3$).

Show that the linear system $|-2K_V-B|$ has no base curves after blowing
up of $B$. It is enough to check that there is an element $G\in |-2K_V|$
that separates the sheets of the double cover. It follows from the
Riemann-Roch theorem for $(-2K_V)$. Indeed, since the canonical divisor
is {\it big} and {\it nef}, we have
$\chi(V,{\cal O}(-2K_V))=H^0(V,{\cal O}(-2K_V))$. Further, $c_2(V)=10s+24f$,
and a direct calculations yeilds
$$
       H^0(V,{\cal O}(-2K_V))=15.
$$
At the same time, we have $H^0(Q,{\cal O}(2L))=14$ for the cone, where
$L$ is the class of a hyperplane section (recall that
$K_V=(\pi\circ\phi_V)^*(-L)$).

Take a general element $G\in |-2K_V-B|$. Let $\psi:V'\to V$ be the blow
up of the curve $B$ with an exceptional divisor $E$. If $\alpha=2$, then
$$
   \psi^{-1}D_1\circ\psi^{-1}D_2\circ\psi^{-1}G=8n^2-4n\nu-4\nu^2 < 0.
$$
For $\alpha=3$ we have
$$
   \psi^{-1}D_1\circ\psi^{-1}D_2\circ\psi^{-1}G=8n^2-6n\nu-5\nu^2 < 0.
$$
But we have a contradiction with
$\psi^{-1}D_1\circ\psi^{-1}D_2\circ\psi^{-1}G\ge 0$.

So, $B$ is a horizontal line. Moreover, assume that $B\ne s_1$ and
$B\ne s_2$. Then $B$ is in a fiber of $\pi_U:U\to\POn$. At the same time,
$B$ is a maximal curve of the linear system $|-nK_U|$ on $U$. But this is
impossible because of {\it(i)}.

\noindent {\it (iii)} Let $l_1$ and $l_2$ be maximal curves, $\nu_1$
and $\nu_2$ corresponding multiplicities of $|D|$ in them.

Almost each fiber $S\in |F_1|$ contains a curve $C_S\in |-K_S|$ that is not
base for $|D|$ and pass through the intersection points of $S$ and those
curves. But we have a contradiction then:
$$
    2n=C_S\circ D\ge\nu_1+\nu_2>2n
$$
The similar arguments prove a case that $l_1$ is infinitely near $l_1$.
The proposition is proved.

Note that this completes the proof of theorem \ref{main_theorem}.
Indeed, let $l$ be a maximal curve. Applying to $V$
the birational automorphisms $\tau_l$, $\tau_1$ or $\tau_2$ (if $l$
coincide with $s_1$ or $s_2$) and lowering the linear system
$|D'|=\tau_l^*|D|$ to $X$, we get $\alpha(X,|D'|) < \alpha(X,|D|)$
as it follows from (\ref{tau_l}) and (\ref{tau_1}).
If $|D|\subset |-nK_X|$, any maximal singularity only lie over one of the
singular points of $X$ because of {\it (ii)}. Thus, we apply
$\tau_1$ or $\tau_2$ in this case.

\section{Proof of corollary \ref{main_corollary}}
{\it (i)} Assume the converse, i.e., $X$ is rational and there is a
birational map
$$
     \chi: X-\!\! -\!\!\to \PT
$$
Theorem \ref{main_theorem} implies that there exists a birational map
$\chi':X\to\PT$ such that
$$
    \alpha(X,(\chi')^{-1}{\cal H}) \le \alpha(\PT,{\cal H})=\frac14,
$$
where $\cal H$ is the complete linear system of hyperplanes in $\PT$.
Then the adjunction treshold $\alpha(X,(\chi')^{-1}{\cal H})=0$, i.e.,
the linear system $(\chi')^{-1}{\cal H}$ is either $|F_1|$ or $|F_2|$.
This is impossible by the dimension reasons.

The similar arguments proves that $X$ is not birationally isomorphic to
a conic bundle:  we can choose ${\cal H}$ is a lift of a sufficiently
ample linear system from a base.

\noindent {\it (ii)} Let we have a birational morphism
$$
     \chi: X-\!\! -\!\!\to X.
$$
Lift it to $V$:
$$
     \chi_V: V-\!\! -\!\!\to V
$$
Let ${\cal H}=|-nK_V+mF_1|$, $n,m > 0$, be a very ample linear system,
for example, $|n(-K_V+F_1)|$, where $n$ is big enough.

There is a birational automorphism $\tilde\chi_V$ such that
the adjunction treshold $n_1=\alpha(X,|D|)$ of the linear system
$|D|=\mu^{-1}{\cal H}$, where $\mu=\chi_V\circ\tilde\chi_V$,
is not greater than $n=\alpha(X,{\cal H})$.

Suppose $|D|=|-n_1K_V+aF_1|$. We use theorem 4.2 of paper \cite{Corti}:
\begin{theorem}[Corti]
(i) $n_1\ge n$;

\noindent (ii) if $K_V+\frac1{n_1}|D|$ is canonical and nef,
then $\mu$ is an isomorphism.
\end{theorem}

So, $n=n_1$ in our case. Further, it is easy to see that $a\ge 0$.
Indeed, assume the converse, then $\dim|D| < \dim{\cal H}$, but this is
impossible.

Thus, $|D|$ is under condition {\it (ii)} of Corti's theorem, whence
$\mu$ is an isomorphism. Lowering to $X$, we see that $\chi\circ\tilde\chi$
is an isomorphism, too.

The proposition in section \ref{MC} proves the uniqueness of the process
of "untwisting" maximal singularities, from which it follows that
there are no any relations in $B(X)$. This implies (\ref{sl1}).

Finally, to prove (\ref{sl2}), note that any regular automorphism of the
double cone is induced by a regular automorphism of $Q$ that preserve
the ramification divisor $R_Q$. In general case there are no such
automorphisms of $Q$. This proves (\ref{sl2}), just note that the groups
$Bir(V_{\eta})$ и $Bir(U_{\eta})$ are completely decribed in \cite{Isk4},
theorem 2.6.

{\large Steklov Mathematical Institute}

{\large Russian Academy of Sciences}

{\large e-mail:  grin @ mi.ras.ru}
\end{document}